\documentclass[12pt,a4paper]{amsart}
\usepackage{amssymb,graphicx}

\newcommand{\re}{\mathbb R}

\begin{document}

\title[Metric Properties of Conflict Sets] {Metric Properties of
Conflict Sets}

\author[L. Birbrair]{Lev Birbrair}
\address {Departamento de Matematica UFC}  \email{birb@ufc.br}

\author[D. Siersma]{Dirk Siersma}
\address{Department of Mathematics University of Utrecht}
\email{siersma@math.uu.nl}

\date{\today}

\markboth {L. Birbrair D. Siersma }{Metric Properties of Conflict
Sets}

\maketitle

\begin {abstract}

 In this paper we show that the tangent cone of a conflict set in $R^n$ is a linear
 affine cone over a conflict set of smaller dimension and has dimension $n-1$. Moreover
 we give an example where the conflict sets is not normally embedded and not
 locally bi-Lipschitz equivalent to the corresponding tangent cone.

\end {abstract}

\section{Introduction}

\bigskip

Singularities of conflict sets of collections of disjoint subsets of
$\re^n$ is one of natural objects of Singularity Theory. Conflict
sets are the boundary of territories of disjoint closed subset with
respect to the nearest distance criterium. The investigation of
conflict sets was initiated by Y.Yomdin \cite{Y}. J.Damon \cite{Da},
P.Giblin and V.M.Zakalyukin \cite{GZ} and others pointed out that
the theory of conflict sets is closely related to other important
objects in Singularity Theory: cut loci, medial axes, wave fronts.
The results of \cite{Da}, \cite{S}, \cite{vM} are devoted to
differential geometry of conflict sets.

Here we study properties of general (not necessary generic)
singularities of conflict set from a metric view-point. Our
restriction is so-called "definability". We suppose that all the
sets appearing in our investigation are definable in some o-minimal
structure \cite{D}. If a reader is not familiar with "o-minimal"
language he can suppose that all the sets are semialgebraic or
subanalytic. If the sets $X_i$ are definable then the same is true
for their conflict sets, Thus, all the "good" topological properties
hold and the Hausdorff limits and tangent cones are well defined.

The main statement of this paper is the structural theorem about the
tangent cone of a conflict set. We prove that it is a linear affine
cone over a conflict set of smaller dimension (Theorem 2.2). As a
corollary of this result we obtain that, for each set $X_i$ of the
collection $\{ X_1,X_2, \ldots, X_k \}$ in $\re^2$, the set $Ter \,
X_i$ does not have the cusp-like regions. Moreover, we show that the
tangent cone to $Ter \, X_i$ has a dimension of the ambient space
and the tangent cone of $Conf\, (X)$ has a dimension $n-1$. A
natural question is the following. Is it true or not that conflict
sets have "metrically conic" structure near a singular point, and
more special, a conflict set is locally bi-Lipschitz equivalent to
the corresponding tangent cone?

The answer is NO: In section 3, we present an example of a
collection of sets $\{ X_1,X_2 \} \subset \re^3$ such that the
conflict set of $X_1$ and $X_2$ is not locally homeomorphic to its
tangent cone and not normally embedded in $\re^3$.

\bigskip

\section{Tangent cones of conflict sets}

\medskip

Let $M$ be a metric space and let $X=\{ X_1, \ldots, X_k \}$ be a
finite collection of closed disjoint nonempty subsets of $M$. We
define a {\it territory \/} of a subset $X_i \in X$ with respect to
the space $M$ and the collection $X$ in the following way:
$$
Ter_M(X_i,X)=\{ x \in M \quad \text{such that} \quad \forall j \quad
d(x,X_i) \le d(x,X_j) \}.
$$
We define a {\it conflict set \/} of the collection $X$ with respect
to $M$ as follows:
$$
Conf_M(X) = \{ x \in Ter_M(X_i,X) \cap Ter_M(X_j,X) \quad \text{for
some} \quad i \neq j \}.
$$

Let $A$ be an o-minimal structure over $\re$. Let $X$ be a
collection of definable subsets of a definable set $Y \subset \re^n$
in an o-minimal structure $A$. Then $Conf_Y(X)$ and $Ter_Y(X_i,X)$
are definable in $A$ sets.

In this paper we are going to suppose that the space $\re^n$ is
equipped with the Euclidean metric.

\bigskip

\noindent{\bf Proposition 2.1.}{\it Let $X = \{ X_1,X_2,X_3, \ldots,
X_k \}$ be a finite collection of closed and definable in an
o-minimal structure $A$ subsets  of $\re^n$. Then $\dim \,
Conf_{\re^n}(X) = n-1$. \/}

\medskip

\noindent PROOF. Let us first show that if $x_0 \in Conf_{\re^n}(X)$
then there exists a number $j$ such that $x_0 \in \partial
Ter_M(X_j,X)$. Let $r_0 = \min \, d(x_0,X_i)$. Let $\widetilde X =
\{ \widetilde X_1, \widetilde X_2, \ldots , \widetilde X_k \}$ be a
collection of sets (called "supports") defined as follows:
$\widetilde X_i = X_i \cap S_{x_0,r_0}$. Since $x_0 \in
Conf_{\re^n}(X)$ we can suppose that there exist two numbers $j_1$
and $j_2$ such that $\widetilde X_{j_1}$ and $\widetilde X_{j_2}$
are nonempty. Let $x_{j_1} \in \widetilde X_{j_1}$ and $x_{j_2} \in
\widetilde X_{j_2}$. Then the half-open segment $[x_{j_1},x_0)$
belongs to $Ter_{\re^n}(X_{j_1},X)$ and does not belong to
$Ter_{\re^n}(X_{j_2},X)$, the segment $(x_{j_2},x_0)$ belongs to
$Ter_{\re^n}(X_{j_2},X)$ and does not belong to
$Ter_{\re^n}(X_{j_1},X)$. Hence, $x_0$ is a boundary point of
$Ter_{\re^n}(X_{j_1},X)$ and of $Ter_{\re^n}(X_{j_2},X)$. This
argument also proves that the sets $Int (Ter_{\re^n}(X_i,X))$ are
disjoint. Since the sets $Ter_{\re^n}(X_i,X)$ are definable in the
o-minimal structure $A$, we have:$\displaystyle Conf_{\re^n}(X)
\subset{\bigcup}_i \partial (Ter_{\re^n}(X_i,X))$. That is why $\dim
Conf_{\re^n}(X) \le n-1$ (see \cite{D}). From the other hand,
$\displaystyle \re^n = {\bigcup}_i Int (Ter_{\re^n}(X_i,X))) \bigcup
Conf_{\re^n}(X)$. Since the sets $Int (Ter_{\re^n}(X_i,X))$ are
disjoint, we obtain that $\dim Conf_{\re^n}(X) \ge n-1$ (see [4]) . This
proves the proposition.  \hfill $\Box$

\bigskip

\noindent{\bf Remark.} A similar statement is true for collections
of definable subsets of $S^n$ and is not true for collections of
definable subsets of $\re^n$ equipped with New York metric.

\bigskip

Let $M$ be a subset of $\re^n$. The cone over $M$ with respect to
$x_0$ (notation: $C_{x_0}M$) is a union of all rays connecting $x_0$
with all the points $y \in M$.

Let $Y \subset \re^n$ be definable in an o-minimal structure $A$. A
tangent cone $T_{x_0}X$ at a point $x_0 \in X$ is the set of all
tangent vectors $\frac{d}{d\gamma}\big|_{t=0}\gamma(t)$ of all
definable in $A$ arcs $\gamma \colon [0, \varepsilon ) \to Y$ such
that $\gamma(0) = x_0$
  ( see \cite{F}, \cite{B}). If $x_0$ is a smooth point of $Y$ then we obtain the definition of the tangent
space.

We are going to use another definition of the tangent cone (see also
\cite{F}). Let $N_{\varepsilon} (Y)$ be a set defined as follows:
$$
N_{\varepsilon}(Y) = \frac{1}{\varepsilon}\big[ Y\bigcap
S_{x_0,\varepsilon}-x_0 \big] + x_0.
$$
We use the notations: $S_{x_0,\varepsilon}$, for the sphere centered
at $x_0$ with the radius $\varepsilon$; $B_{x_0,\varepsilon}$, for
the closed ball.

Then the following statement is true:
$$
T_{x_0}Y= C_{x_0}({\lim_{{Hausdorff}\varepsilon \to 0}}
N_\varepsilon(Y)).
$$

\bigskip

Let $X = \{ X_1, \ldots, X_k \}$ be a family of closed and disjoint
sets on $S^{n-1}$. Then there is the following relation:
$$
Conf_{\re^n}(X) = C_{o\in \re^n}Conf_{S^{n-1}}(X).
$$
Here we use the standard geodesic metric on $S^{n-1}$. The proof of
this statement is straightforward.

\bigskip

The main result of this paper is the following statement.

\bigskip

\noindent{\bf Theorem 2.2.}{\it Let $X = \{ X_1, \ldots, X_k \}$ be
a collection of definable in $A$ closed subsets of $\re^n$ such that
$X_i \cap X_j = \emptyset$, for $i \neq j$. Let $x_0 \in
Conf_{\re^n}(X)$ and let $r_0 = \underset{i}\to{\min}\, d(x_0,X_i)$.
Let $\widetilde X = \{ \widetilde X_1, \widetilde X_2, \ldots,
\widetilde X_k \}$ be a collection of sets ( called "the supports" )
defined as follows: $\widetilde X_i = X_i \cap S_{x_0,r_0}$. Then
the following identities hold:}

\begin{enumerate}
\item $T_{x_0}(Ter_{\re^n}(X_i,X))=C_{x_0}(Ter_{S_{x_0,r_0}}
(\widetilde X_i, \widetilde X))$.
\item $T_{x_0}(Conf_{\re^n}(X))=C_{x_0}(Conf_{S_{x_0,r_0}}
(\widetilde X))$.
\end{enumerate}

\noindent PROOF. Observe that the statement $2$ follows immediately
from the statement 1 by the definition of conflict sets. Now we are
going to show that the germs of $Ter_{\re^n}(X_i,X)$ and of
$Conf_{\re^n}(X)$ at $x_0$ do not change if we cut the sets $X_i$ by
balls of the radius bigger than $r_0$ centered at $x_0$. Namely, we
prove the following statement:

\medskip

\noindent{\bf Lemma 2.3.} {\it Let $X^{\varepsilon} = \{
X_1^{\varepsilon}, X_2^{\varepsilon}, \ldots, X_k^{\varepsilon} \}$
be a collection of the sets defined as follows: $X_i^{\varepsilon}
=X_i \cap \bar B_{x_0,r_0+\varepsilon}$. Then, for all $i$, the germ
of $Ter_{\re^n}(X_i,X)$ at $x_0$ is the same as the germ of the set
$Ter_{\re^n}(X_i^{\varepsilon},X^{\varepsilon})$ at $x_0$.\/}

\medskip

\noindent PROOF.  Let $z$ be a point of $X_j$ such that
$d(x_0,z)=d(x_0,X_j)$. Take $\delta = \varepsilon /3$. Let $x \in
B_{x_0,\delta}$. Let $y$ be a point of $X_j$ such that $d(x,y)=
d(x,X_j)$. We are going to show that $y$ actually belong to
$X_j^{\varepsilon}$.

First: $d(x,y)=d(x,X_j) \le d(x,z) \le d(x,x_0)+d(x_0,z) \le \delta
+r_0$. Second: suppose that $y \notin X_j^{\varepsilon}$. Then
$d(x,y) > d(x_0,y)-d(x,x_0)=r_0+\varepsilon -\delta$. This is a
contradiction. So, on $B_{x_0,\delta}$ one have:
$d(x,X_j)=d(x,X_j^{\varepsilon})$ and, therefore,
$Ter_{\re^n}(X_i,X) \cap
B_{x_0,\delta}=Ter_{\re^n}(X_i^{\varepsilon},X^{\varepsilon}) \cap
B_{x_0,\delta}$. \hfill $\Box$

We consider "polar coordinates" $(\rho , \phi )$ near the point
$x_0$ defined as follows. Let $x \in \re^n$ be a point. Set $\phi
(x) = \frac{x-x_0}{||x-x_0||}$ and $\rho (x) = d(x,x_0)$.

\bigskip

\noindent{\bf Lemma 2.4. (Shadow lemma).}  {\it
Let $X = \{ X_1, \ldots, X_n \}$ be a family of definable sets such
that, for all $j$, $\phi (X_j) = \phi (\widetilde X_j)$ (we say that
$X$ lies in the shadow of the support of $X$ at $x_0$). Then the
germs at $x_0$ of the sets $Ter_{\re^n}(X_j,X)$ and
$Ter_{\re^n}(\widetilde X_j,\widetilde X)$ are equal. \/}

\medskip

\noindent{PROOF.} Let $x_1 \in B_{x_0,r_0/3}$. Then $d(x_1,X_j) =
d(x_1, \widetilde X_j)$. That is why
$$
Ter_{\re^n}(\widetilde X_j, \widetilde X) \cap B_{x_0,r_0/3} =
Ter_{\re^n}(X_j, X) \cap B_{x_0,r_0/3}.
$$
\hfill $\Box$

\medskip

\noindent END OF THE PROOF OF THEOREM 2.2. Using polar coordinates
we define sets $Y_i^{\varepsilon}$ in the following way: $y \in
Y_i^{\varepsilon}$ if, and only if, $y\in B_{x_0,r_0+\varepsilon},
\quad  r_0\le \rho(y) \le r_0 + \varepsilon$ and there exists $x \in
X_i^{\varepsilon} \quad \text{s.t.} \quad \phi(x) = \phi (y)$.
Clearly, $X_i^{\varepsilon} \subset Y_i^{\varepsilon}$. Let
$Y^{\varepsilon} = \{ Y_1^{\varepsilon},Y_2^{\varepsilon}, \ldots,
Y_k^{\varepsilon} \}$. Let $Z^{\varepsilon} = \{
Z_1^{\varepsilon},Z_2^{\varepsilon}, \ldots, Z_k^{\varepsilon} \}$
be a collection of sets defined as follows: $Z_i^{\varepsilon} =
Y_i^{\varepsilon} \cap S_{x_o,r_0}$. The Hausdorff limits for
$\varepsilon$ tending to zero of the families $Y_i^{\varepsilon}$
and $Z_i^{\varepsilon}$ are equal to $\widetilde X_i$. If $x \in
B_{x_o, r_0/3}$ and $d(x,Y_i^{\varepsilon}) = d(x,y)$, for some $y
\in Y_i^{\varepsilon}$, then $y \in Z_i^{\varepsilon}$. Clearly, for
small $\varepsilon$, we have: $Y_i^{\varepsilon} \cap
Y_j^{\varepsilon} = \emptyset$, for $i \neq j$, and, thus,
$Z_i^{\varepsilon} \cap Z_j^{\varepsilon} = \emptyset$. Since all
the sets $Z_i^{\varepsilon}$ belong to $S_{x_0,r_0}$, then
$Ter_{\re^n}(Z_i^{\varepsilon},Z^{\varepsilon})=
C_{x_0}Ter_{S_{x_0,r_0}}(Z_i^{\varepsilon},Z^{\varepsilon})$.
Moreover, the germ of the set $Ter_{\re^n}(Y_i^{\varepsilon},
Y^{\varepsilon})$ at the point $x_0$ is the same as the germ of
$Ter_{\re^n}(Z_i^{\varepsilon}, Z^{\varepsilon})$ [Shadow Lemma].

Let us consider the collection of sets $W^{\varepsilon} = \{
W_1^{\varepsilon},W_2^{\varepsilon}, \ldots, W_k^{\varepsilon} \}$
and $V^{\varepsilon} = \{ V_1^{\varepsilon},V_2^{\varepsilon},
\ldots, V_k^{\varepsilon} \}$ where $W_j^{\varepsilon} =
Y_j^{\varepsilon}$, for $i \neq j$, and $W_i^{\varepsilon}=
\widetilde X_i$, $V_j^{\varepsilon}=\widetilde X_j$, for  $i \neq
j$, and $V_i^{\varepsilon}=Y_i^{\varepsilon}$. Since $\widetilde X_s
\subset X_s^{\varepsilon} \subset Y_s^{\varepsilon}$, for any $s$,
we obtain the following inclusions:
$$
Ter_{\re^n}(X_i^{\varepsilon},X^{\varepsilon}) \subset
Ter_{\re^n}(Y_i^{\varepsilon},V^{\varepsilon}),
$$
 $\quad\quad \quad \quad \quad \quad \quad\quad\quad \quad \quad \quad \quad \quad\quad \quad \quad \quad
  \quad\quad\quad \quad \quad \quad \quad \quad (1)$
$$
Ter_{\re^n}(\widetilde X_i,W^{\varepsilon}) \subset
Ter_{\re^n}(X_i^{\varepsilon},X^{\varepsilon}).
$$

Let $\widetilde W^{\varepsilon}$ and $\widetilde V^{\varepsilon}$ be
the collections of sets defined as follows: $\widetilde
W^{\varepsilon} =  W^{\varepsilon} \cap S_{x_0,r_0}$ and $\widetilde
V^{\varepsilon} =  V^{\varepsilon} \cap S_{x_0,r_0}$. Note, that the
germs of the sets $Ter_{\re^n}(V_j^{\varepsilon},V^{\varepsilon})$
and $Ter_{\re^n}(\widetilde V_j^{\varepsilon},\widetilde
V^{\varepsilon})$ at $x_0$ are equal and the germs of
$Ter_{\re^n}(W_j^{\varepsilon},W^{\varepsilon})$ and
$Ter_{\re^n}(\widetilde W_j^{\varepsilon},\widetilde
W^{\varepsilon})$ are also equal (Shadow lemma). The sets
$Ter_{\re^n}(\widetilde V_i^{\varepsilon},\widetilde
V^{\varepsilon})$ and $Ter_{\re^n}(\widetilde
W^{\varepsilon},\widetilde W)$ are purely conic., i.e.
$$
Ter_{\re^n}(\widetilde V_i^{\varepsilon},\widetilde V^{\varepsilon})
= C_{x_0}Ter_{S_{x_0,r_0}}(\widetilde V_i^{\varepsilon},\widetilde
V^{\varepsilon})
$$
and
$$
Ter_{\re^n}(\widetilde W_i^{\varepsilon},\widetilde W^{\varepsilon})
= C_{x_0,r_0}Ter_{S_{x_0,r_0}}(\widetilde
W_i^{\varepsilon},\widetilde W^{\varepsilon}).
$$
Hence,
$$
T_{x_0}Ter_{\re^n}(\widetilde V_i^{\varepsilon},\widetilde
V^{\varepsilon})=Ter_{\re^n}(\widetilde V_i^{\varepsilon},\widetilde
V^{\varepsilon})
$$
and
$$
T_{x_0}Ter_{\re^n}(\widetilde W_i^{\varepsilon},\widetilde
W^{\varepsilon})=Ter_{\re^n}(\widetilde W_i^{\varepsilon},\widetilde
W^{\varepsilon}).
$$

Using the inclusions $(1)$ we obtain:

$$
T_{x_0}(Ter_{\re^n}(X_i,X))=
T_{x_0}(Ter_{\re^n}(X_i^{\varepsilon},X^{\varepsilon})) \subset
Ter_{\re^n}(\widetilde V_i^{\varepsilon},\widetilde
V^{\varepsilon}),\quad \quad (2)
$$
$$
Ter_{\re^n}(\widetilde W_i^{\varepsilon},\widetilde W^{\varepsilon})
\subset
T_{x_0}(Ter_{\re^n}(X_i^{\varepsilon},X^{\varepsilon}))=T_{x_0}(Ter_{\re^n}(X_i,X)).
$$

Taking the Hausdorff limit in the inclusions $(2)$ we obtain
$$
T_{x_0}(Ter_{\re^n}(X_i,X)) \subset
C_{x_0}(Ter_{S_{x_0,r_0}}(\widetilde X_i,\widetilde X))
$$
and
$$
C_{x_0}(Ter_{S_{x_0,r_0}}(\widetilde X_i,\widetilde X)) \subset
T_{x_0}(Ter_{\re^n}(X_i,X)).
$$
This proves the theorem. \hfill $\Box$

\bigskip

\noindent{\bf Proposition 2.6.} {\it Let $X = \{ X_1, \ldots, X_k
\}$ be a collection of definable in $A$ sets in $\re^n$.  Then}

\begin{enumerate}
\item {\it $T_{x_0}(Ter_{\re^n}(X_i,X))$ has a nonempty interior, for all $i$ such that
$x_0 \in Ter (X_i,X)$.}
\item {\it If $x_0 \in Conf_{\re^n}(X)$ then
$\dim (T_{x_0}Conf_{\re^n}(X))=n-1$.}
\end{enumerate}

\noindent PROOF. [1] If $\widetilde X_i \neq \emptyset$ then
$Ter_{S^{n-1}}(\widetilde X_i,\widetilde X)$ has a nonempty
interior. Thus, by Theorem 2.2, $T_{x_0}(Ter_{\re^n}(X_i,X))$ has a
nonempty interior. [2] By Proposition 2.1 (see also the remark),
$\dim Conf_{S^{n-1}}(\widetilde X)=n-2$. Hence, $\dim
(T_{x_0}Conf_{\re^n}(X))=n-1$. \hfill $\Box$

\bigskip

\noindent{\bf Theorem 2.7 ("no cusp" property in $\re^2$).} {\it Let
$X = \{ X_1, \ldots, X_k \}$ be a collection of definable in $A$
sets on $\re^2$. Let $Y=Conf_{\re^2}(X)$. Let $y_0 \in Y$. Then the
germ of $Y$ at $y_0$ is a collection of definable in $A$ arcs $\{
\gamma_1, \gamma_2, \ldots, \gamma_s \}$ such that $y_0$ belongs to
each $\gamma_i$ and the unit tangent vectors of $\gamma_i$ and
$\gamma_j$ are different, for $i \neq j$. \/}

\medskip

\noindent PROOF. Let $\gamma_1 \colon [0,\varepsilon ) \to
Conf_{\re^2}(X)$ be a definable nonconstant arc such that $\gamma_1
(0) = y_0$ and $| \gamma_1(t) - y_0 |=t$. Since the sets $X_i$ are
definable in the o-minimal structure $A$ we can find two sets
$X_1,X_2 \in X$ such that $\gamma_1 \subset Conf_{\re^2}(\bar X)$
where $\bar X = \{ X_1,X_2 \}$. Let $l$ be the tangent ray to
$\gamma_1$ at $y_0$. By Theorem 2.2, $l$ is a bisector ray of the
angle defined by points $x_1 \in \widetilde X_1, \quad x_2 \in
\widetilde X_2$ and $y_0$, where $x_1, x_2$ are boundary points on
$S_{y_0,r_0}$ of the supporting sets $\widetilde X_1$ and
$\widetilde X_2$ on $S_{y_0,r_0}$.

Let $\delta > 0$ be a sufficiently small number such that the sets
$X_1^{\delta}$ and $X_2^{\delta}$ - their radial projections to the
supporting circle are disjoint. Let $\gamma_2 \colon [0,\varepsilon
) \to Conf_{\re^2}(\bar X)$ be another definable in $A$ arc such
that $|\gamma_2(t)-y_0 |=t$ and the germs at $y_0$ of the sets
$\Gamma_1 = \gamma_1([0,\varepsilon))$ and $\Gamma_2 =
\gamma_2([0,\varepsilon))$ are different and $\Gamma_2$ is also
tangent to $l$ at $y_0$. By Arc Selection Lemma (see \cite{D}),
there exist two pairs of arcs $\alpha_1,\alpha_2 \colon
[0,\varepsilon) \to X_1^{\delta}$ and $\beta_1,\beta_2 \colon
[0,\varepsilon) \to X_2^{\delta}$ such that
$d(\gamma_1(t),X_1)=|\gamma_1(t)-\alpha_1(t)|, \quad
d(\gamma_2(t),X_1)=|\gamma_2(t)-\alpha_2(t)|, \quad
d(\gamma_1(t),X_2)=|\gamma_1(t)-\beta_1(t)|$ and
$d(\gamma_2(t),X_2)=|\gamma_2(t)-\beta_2(t)|$. If $\alpha_1(t) =
\alpha_2(t)=x_1$ and $\beta_1(t)=\beta_2(t)=x_2$, then the germs of
$\gamma_1$ and $\gamma_2$ at $y_0$ are equal to the germ of $l$ at
$y_0$. Thus, we can suppose that $\alpha_1(t) \neq Const$, for small
$t$, and that, for small $t \neq 0$ we have $d(\gamma_1(t),X_1) >
d(\gamma_2(t),X_1)$.

Take $t>0$ sufficiently small. The segment connecting $\gamma_1(t)$
and $\beta_1(t)$ intersects the arc $\gamma_2(t)$. Let $z$ be an
intersection point. Observe that $\beta_1(t)$ realizes the shortest
distance between $z$ and $X_1$.  Since $t$ is small and $\beta_1$ is
a definable  in $A$ arc, then we can suppose that $\beta_1(t)$ is a
smooth point.

Consider now the circle $S_{\gamma_1(t),|\gamma_1(t)-\beta_1(t)|}$
and the circle $S_{z,|z-\beta_1(t)|}$. These circles are tangent at
the point $\beta_1(t)$. That is why the ball with the center at $z$
and the radius $|z-\beta_1(t)|$ does not contain any point of $X_1$.
But it means that $z$ does not belong to $Conf_{\re^2}(\bar X)$. It
is a contradiction. \hfill $\Box$

\bigskip

\section{An example of not normally embedded conflict set}

\medskip

Here we are going to construct an example of a family of sets
$X_1,X_2 \in \re^3$ satisfying the following conditions:
\begin{enumerate}
\item[a)] $Conf_{\re^3} (\{ X_1,X_2 \} )$ is not normally embedded in $\re^3$.
\item[b)] There exists a point $x_0 \in Conf_{\re^3}(\{ X_1,X_2 \})$
such that, for small $r$, we have that $B_{x_0,r} \cap
Conf_{\re^3}(\{ X_1,X_2\})$ is not homeomorphic to
$T_{x_0}(Conf_{\re^3}(\{ X_1,X_2\} ))$.
\end{enumerate}

\noindent{\bf Example.} Consider the space $\re^3$ with coordinates
$(x_1,x_2,x_3)$. Let $X_1 \subset \re^3$ be a union of the
hyperplanes: $x_3=1$ and $x_3=-1$. Let $X_2$ be a union of the
points: $a_1=(1,0,0)$ and $a_2=(-1,0,0)$.

\medskip

\noindent{\bf Theorem 3.1.}{\it The set $Conf_{\re^3}(\{ X_1,X_2\}
)$ satisfies the conditions $a)$ and $b)$ described above. \/}

\medskip

\noindent PROOF. The set $Conf_{\re^3}(\{ X_1,X_2\} )$ can be
obtained as follows. Let $Y_1 \subset \re^2$ be the conflict set af
the point $a_1 = (1,0)$ and the union of straight lines $x_3=1$ and
$x_3=-1$. Observe, that here we consider $\re^2$ with coordinates
$(x_1,x_3)$. The set $Y_1$ is a union of a part of the parabola
defined by the point $a_1$ and the line $x_3=1$ situated above the
line $x_3=0$ and a part of the parabola defined by the same point
$a_1$ and the line $x_3=-1$ situated below the line $x_3=0$.

Let $Y_2 \subset \re^2$ be the conflict set of the point $a_2 =
(-1,0)$ and the union of the lines $x_3=1$ and $x_3=-1$. The set
$Y_2$ can be obtained from $Y_1$ by the transformation: $(x_1,x_3)
\to (-x_1,x_3)$. The set $Conf_{\re^3}(\{ X_1,X_2\} )$ can be
obtained as a union of the revolution surface of $Y_1$ with respect
to the straight line $x_1=1,\quad x_2=0$ and the revolution surface
of $Y_2$ with respect to the line $x_1=-1, \quad x_2=0$. The
intersection of the set $Conf_{\re^3}(\{ X_1,X_2\} )$ with the plane
$x_3=0$ is a union of two metric copies of $S^1$. These circles are
tangent at the origin. That is why the set $Conf_{\re^3}(\{
X_1,X_2\} )$ is not normally embedded. The tangent cone
$T_0Conf_{\re^3}(\{ X_1,X_2\} )$ is a union of two planes
intersecting transversally. The germ of the set $Conf_{\re^3}(\{
X_1,X_2\} )$ is homeomorphic to the quotient space of the disjoint
union of two copies of $\re^2$ by the identification of the two
origins. \hfill $\Box$

\bigskip

\noindent{Asknowledgements.} The first author was supported by CNPq
grant N 300985/93-2, the second author was supported by CNPq grant N
420108/2005-8.

We would like to thank Y.Yomdin, T.Mostowski, A.Fernandes and M. van
Manen for their interest on this work.

\end{document}